\overfullrule=0pt
\centerline {\bf Multiple critical points in closed sets via minimax theorems}\par
\bigskip
\bigskip
\centerline {BIAGIO RICCERI}\par
\bigskip
\bigskip
{\bf Abstract.} In this paper, we apply our minimax theory ([4], [5], [6]) with the one developed by A. Moameni in [2] to formalize a general
scheme giving the multiplicity of critical points. Here is a sample of application of the scheme to a critical elliptic problem:
\medskip
THEOREM. - {\it Let $\Omega\subset {\bf R}^n$ ($n\geq 3$) be
a smooth bounded domain and let $1<q<2\leq p<{{2n}\over {n-2}}$.\par
Then, for every $r, \nu>0$, there exists $\lambda^*>0$ with the following property: for every $\lambda\in ]0,\lambda^*[$, $\mu\in ]-\lambda^*,\lambda^*[$, and
for every convex dense set $S\subset H^{-1}(\Omega)$, there exists $\tilde\varphi\in S$, with $\|\tilde\varphi\|_{H^{-1}(\Omega)}<r$, such that
the problem
$$\cases{-\Delta u=\lambda(|u|^{{{4}\over {n-2}}}u+\nu |u|^{q-2}u+\mu|u|^{p-2}u+\tilde\varphi) & in $\Omega$\cr & \cr
u=0 & on $\partial\Omega$\cr}$$
has at least two solutions whose norms in $H^1_0(\Omega)$ are less than or equal to $r$.}\par
\bigskip
\bigskip
{\bf 2020 Mathematics Subjet Classification.} 35B38; 35J20; 49J35.
\bigskip
\bigskip
\bigskip
\bigskip
{\bf 1. Introduction}\par
\bigskip
Throughout the sequel, $E$ is a real Banach space and $X\subset E$ is a closed set.\par
\smallskip
Let $g:E\to {\bf R}$ be a functional which is G\^ateaux differentiable at
each point of $X$. A basic question is: under which additional assumptions a global minimum $x_0$ of $g_{|X}$ is also a critical point of $g$
(that is, $g'(x_0)=0$) ? \par
\smallskip
Quite recently, a very remarkable contribution to such a topic was provided by A. Moameni in [2].\par
\smallskip
 We now recall
two results from [2] in a simplified form which is enough for our purposes.\par
\smallskip
 We first recall the following definition. If $\Phi, \Psi:E\to
{\bf R}$ are G\^ateaux differentiable at each point of $X$, we say that the triple $(X,\Phi,\Psi)$ satisfies the pointwise invariant
condition at a point $u_0\in X$ if at least one of the following two facts holds:\par
\smallskip
\noindent
-\hskip 5pt there exist $v_0\in X$ and a lower semicontinuous convex G\^ateaux differentiable functional
$\gamma:E\to {\bf R}$ such that
$$\Phi'(u_0)+\gamma'(u_0)=\Psi'(v_0)+\gamma'(v_0)\ ;$$
-\hskip 5pt there exist $v_0\in X$ and a lower semicontinuous convex G\^ateaux differentiable functional
$\gamma:E\to {\bf R}$, with $\gamma'(0)=0$, such that
$$\Phi'(u_0)=\Psi'(v_0)+\gamma'(v_0-u_0)\ .$$
\medskip
THEOREM 1.A. - {\it Let $X$ be also convex. Let $\Phi:E\to {\bf R}$ be a $C^1$ functional and let $\Psi:E\to {\bf R}$ be
a lower semicontinuous convex functional which is G\^ateaux differentiable at each point of $X$. Assume that $u_0\in X$
is a global minimum of $(\Psi-\Phi)_{|X}$ such that the triple $(X,\Phi,\Psi)$ satisfies the pointwise invariant
condition at the point $u_0$.\par
Then, $u_0$ is a critical point of $\Psi-\Phi$.}\par
\medskip
THEOREM 1.B. - {\it Let $\Phi, \Psi:E\to {\bf R}$ be two lower semicontinuous convex functionals which are G\^ateaux differentiable
at each point of $X$ and let $u_0\in X$ be a global minimum of $(\Psi-\Phi)_{|X}$ such that the triple $(X,\Phi,\Psi)$ satisfies the pointwise invariant
condition at the point $u_0$.\par
Then, $u_0$ is a critical point of $\Psi-\Phi$.}\par
\medskip
On the other hand, in [4], [5], [6], we established a general result which, for a suitable family of real-valued functions, ensures the existence of
a function of the family which possesses at least two global minima. Here is the statement:\par
\medskip
THEOREM 1.C. - {\it Let $V$ be a topological space, let $Y$ be a convex set in a topological vector space and let $f:V\times Y\to {\bf R}$
 continuous and quasi-concave in $Y$. Also, assume that
$$\sup_Y\inf_Vf<\inf_V\sup_Yf\ .$$
Moreover, let $\varphi:V\to {\bf R}$ be such that
$$\sup_V\varphi-\inf_V\varphi<\inf_V\sup_Yf-\sup_Y\inf_Vf$$
and that $f(\cdot,y)+\varphi(\cdot)$ is lower semicontinuous and inf-compact for each $y\in Y$.\par
Then, for each convex set $S\subseteq Y$, dense in $Y$, there exists $\tilde y\in S$ such that the function $f(\cdot,\tilde y)+\varphi(\cdot)$ has at least two global minima.}\par
\medskip
The aim of this short note is simply to present an initial set of applications of the general scheme which arises from
combining Theorems 1.A, 1.B and 1.C together, as formalized by Theorems 1.1 and 1.2 below.\par
\medskip
THEOREM 1.1. - {\it Let $X$ be also convex and let $Y$ be a convex set in a topological vector space. Let $I, J: E\times Y\to {\bf R}$ and
$\alpha, \beta:E\to {\bf R}$. Assume that $\alpha$ and $I(\cdot,y)$ (for each $y\in Y$) are G\^ateaux differentiable at each point of
$X$, and that $\beta$ and $J(\cdot,y)$ (for each $y\in Y$) are $C^1$ in $E$. Assume also that
$$\sup_Y\inf_X(I-J)<\inf_X\sup_Y(I-J)$$
and that
$$\sup_X(\alpha+\beta)-\inf_X(\alpha+\beta)<\inf_X\sup_Y(I-J)-\sup_Y\inf_X(I-J)\ .$$
Moreover, suppose the following other conditions are satisfied:\par
\noindent
$(i_1)$\hskip 5pt for each $y\in Y$, the function $I(\cdot,y)+\alpha(\cdot)$ is lower semicontinuous and convex in $E$\ ;\par
\noindent
$(i_2)$\hskip 5pt there exists a topology $\tau$ on $X$ such that, for each $y\in Y$, the function
$I(\cdot,y)-J(\cdot,y)+\alpha(\cdot)+\beta(\cdot)$ is $\tau$-lower semicontinuous and $\tau$-inf-compact in $X$\ ;\par
\noindent
$(i_3)$\hskip 5pt for each $x\in X$, the function $I(x,\cdot)-J(x,\cdot)$ is continuous and quasi-concave in $Y$ ;\par
\noindent
$(i_4)$\hskip 5pt for each $y\in Y$, and for each $u\in X$ such that
$$I(u,y)-J(u,y)+\alpha(u)+\beta(u)=\inf_{x\in X}(I(x,y)-J(x,y)+\alpha(x)+\beta(x))$$
the triple $(X, J(\cdot,y)-\beta(\cdot),I(\cdot,y)+\alpha(\cdot))$ satisfies the pointwise invariance condition at the point $u$.\par
Then, for each convex set $S\subseteq Y$ , dense in $Y$, there exists $\tilde y\in S$ such that the equation
$$I'_x(x,\tilde y)-J'_x(x,\tilde y)+\alpha'(x)+\beta'(x)=0$$
has at least two solutions in $X$.}\par
\medskip
THEOREM 1.2. - {\it Let $Y$ be a convex set in a topological vector space. Let $I, J: E\times Y\to {\bf R}$ and
$\alpha, \beta:X\to {\bf R}$. Assume that $\alpha, \beta$, $I(\cdot,y)$ and $J(\cdot,y)$ (for each $y\in Y$) are G\^ateaux differentiable at each point of $X$ Assume also that
$$\sup_Y\inf_X(I-J)<\inf_X\sup_Y(I-J)$$
and that
$$\sup_X(\alpha+\beta)-\inf_X(\alpha+\beta)<\inf_X\sup_Y(I-J)-\sup_Y\inf_X(I-J)\ .$$
Moreover, besides conditions $(i_1)-(i_4)$ of Theorem 1.1, assume that, for each $y\in Y$, the function $J(\cdot,y)-\beta(\cdot)$ is lower semicontinuous and convex in $E$.\par
Then, the conclusion of Theorem 1.1 holds.}\par
\medskip
On the basis of Theorems 1.A, 1.B and 1.C the proofs of Theorem 1.1 and 1.2 are immediate. In any case, we formalize the proof of
Theorem 1.1. That of Theorem 1.2 is the same, with the only difference of using Theorem 1.B instead of Theorem 1.A.
\medskip
{\it Proof of Theorem 1.1.} Consider the functions $f:X\times Y\to {\bf R}$, $\varphi:X\to {\bf R}$ defined by
$$f(x,y)=I(x,y)-J(x,y)\ ,$$
$$\varphi(x)=\alpha(x)+\beta(x)$$
for all $(x,y)\in X\times Y$. Taking $V=(X,\tau)$, we then see that the assumptions of Theorem 1.C are satisfied. Consequently, for
each convex set $S\subseteq Y$, dense in $Y$, there exists $\tilde y\in S$ such that $f(\cdot,\tilde y)+\varphi(\cdot)$ has at least
two global minima in $X$, say $x_1, x_2$. Now, we consider the functions $\Phi, \Psi:E\to {\bf R}$ defined by
$$\Phi(x)=J(x,\tilde y)-\beta(x)\ ,$$
$$\Psi(x)=I(x,\tilde y)+\alpha(x)$$
for all $x\in E$. So, $\Phi$ is $C^1$ and $\Psi$ is lower semicontinuous, convex and G\^ateaux differentiable at each point of $X$. Since
$$\Psi(x)-\Phi(x)=f(x,\tilde y)+\varphi(x)\ ,$$
the triple $(X,\Phi,\Psi)$ satisfies the pointwise invariant property at $x_1$ and $x_2$. Consequently, Theorem 1.A ensures that $x_1, x_2$
are critical points of $\Psi-\Phi$, and the proof is complete.\hfill $\bigtriangleup$
\bigskip
{\bf 2. Applications}\par
\bigskip
We first recall a well-known proposition.\par
\medskip
PROPOSITION 2.A. - {\it Let $V, Y$ be two topological spaces and let $f:V\times Y\to {\bf R}$ be lower semicontinuous in $V$ and
upper semicontinuous in $Y$. Moreover, assume that there exist $x_0\in V$ and $y_0\in Y$ so that $f(\cdot,y_0)$ is inf-compact and
$f(x_0,\cdot)$ is sup-compact.\par
Then, the following assertions are equivalent:\par
\noindent 
$(c_1)$\hskip 5pt $\sup_Y\inf_Vf=\inf_V\sup_Yf$\ ;\par
\noindent
$(c_2)$\hskip 5pt there exists $(x^*,y^*)\in V\times Y$ such that
$$f(x^*,y^*)=\inf_{x\in V}f(x,y^*)=\sup_{y\in Y}f(x^*,y)\ .$$}
\medskip
For each $r>0$, we set
$$B_r=\{x\in E : \|x\|\leq r\}\ .$$
\medskip
THEOREM 2.1. - {\it Let $E$ be a Hilbert space, let $r>0$ and let $P:E\to {\bf R}$ be a $C^1$ 
even functional, with $P(0)=0$, such that
$$\limsup_{x\to 0}{{P(x)}\over {\|x\|^2}}=+\infty \eqno{(2.1)}$$
and the restriction to $B_r$ of the functional $x\to \sigma\|x\|^2-P(x)$ is sequentially weakly lower semicontinuous for some 
$\sigma>0$.\par
Then, for each $M\geq 0$ and for each $\lambda$ satisfying
$$0<\lambda\leq \min\left\{{{r}\over {\sup_{x\in B_r}\|x+P'(x)\|+M}}, {{1}\over {2\sigma}}\right\}\ ,$$
there exists $\delta^*>0$ with the following property: 
for each convex set $S\subseteq B_r$, dense in $B_r$, and for each
$C^1$ functional $\gamma:E\to {\bf R}$ satisfying 
$$\sup_{B_r}\gamma-\inf_{B_r}\gamma<\delta^*\ ,$$
$$\sup_{x\in B_r}\|\gamma'(x)\|\leq M\ ,$$
such that the restriction to $B_r$ of the functional $x\to \sigma\|x\|^2-P(x)-\gamma(x)$ is sequentially weakly lower semicontinuous,
there exists $\tilde y\in S$, with $\|\tilde y\|<r$, such that the equation
$$x=\lambda(P'(x)+\gamma'(x)+\tilde y)$$
has at least two solutions lying in $B_r$.}\par
\smallskip
PROOF. Fix $M\geq0$ and fix $\lambda$ so that
$$0<\lambda\leq \min\left\{{{r}\over {\sup_{x\in B_r}\|x+P'(x)\|+M}}, {{1}\over {2\sigma}}\right\}\ .\eqno{(2.2)}$$
Set
$$\rho={{\lambda+1}\over {\lambda}}\ .$$
Now, we claim that
$$\sup_{y\in B_r}\inf_{x\in B_r}\left({{\rho}\over {2}}\|x\|^2-{{1}\over {2}}\|x-y\|^2-P(x)\right)<
\inf_{x\in B_r}\sup_{y\in B_r}\left({{\rho}\over {2}}\|x\|^2-{{1}\over {2}}\|x-y\|^2-P(x)\right)\ .\eqno{(2.3)}$$
Arguing by contradiction, assume the contrary, that is
$$\sup_{y\in B_r}\inf_{x\in B_r}\left({{\rho}\over {2}}\|x\|^2-{{1}\over {2}}\|x-y\|^2-P(x)\right)=
\inf_{x\in B_r}\sup_{y\in B_r}\left({{\rho}\over {2}}\|x\|^2-{{1}\over {2}}\|x-y\|^2-P(x)\right)\ .\eqno{(2.4)}$$
We now apply Proposition 2.A taking $V=Y=B_r$, $B_r$ being endowed with the weak topology, and
$$f(x,y)={{\rho}\over {2}}\|x\|^2-{{1}\over {2}}\|x-y\|^2-P(x)$$
for all $(x,y)\in V\times Y$. The assumptions of Proposition 2.A are satisfied. Indeed, $V, Y$ are weakly compact, $f(x,\cdot)$
is weakly upper semicontinuous in $Y$ for all $x\in V$ and $f(\cdot,y)$ is weakly lower semicontinuous in $X$ for
all $y\in Y$, since, as $\lambda\leq {{1}\over {2\sigma}}$, we have ${{\rho-1}\over {2}}\geq\sigma$.
Consequently, in view of $(2.4)$, there would exist $(x^*,y^*)\in B_r\times B_r$ such that
$${{\rho}\over {2}}\|x^*\|^2-{{1}\over {2}}\|x^*-y^*\|^2-P(x^*)=\inf_{x\in B_r}\left({{\rho}\over {2}}\|x\|^2-{{1}\over {2}}\|x-y^*\|^2-P(x)\right)=\sup_{y\in B_r}\left({{\rho}\over {2}}\|x^*\|^2-{{1}\over {2}}\|x^*-y\|^2-P(x^*)\right)\ .$$
From this, in particular, we obtain
$$-{{1}\over {2}}\|x^*-y^*\|=\sup_{y\in B_r}-{{1}\over {2}}\|x^*-y\|^2=-{{1}\over {2}}\inf_{y\in B_r}\|x^*-y\|^2=0$$
and hence $x^*=y^*$. Then, we have
$${{\rho}\over {2}}\|x^*\|^2-P(x^*)=\inf_{x\in B_r}\left({{\rho}\over {2}}\|x\|^2-{{1}\over {2}}\|x-x^*\|^2-P(x)\right)\ .$$
Notice that $x^*=0$. Indeed, if $x^*\neq 0$,  remembering that $P$ is even, we would have
$${{\rho}\over {2}}\|-x^*\|^2-{{1}\over {2}}\|-2x^*\|^2-P(-x^*)<{{\rho}\over {2}}\|x^*\|^2-P(x^*)\ ,$$
a contradiction. But then, since $P(0)=0$,  we have
$${{\rho}\over {2}}\|x\|^2-{{1}\over {2}}\|x\|^2-P(x)\geq 0$$
for all $x\in B_r$, that is
$${{P(x)}\over {\|x\|^2}}\leq {{1}\over {2}}(\rho-1)$$
for all $x\in B_r\setminus \{0\}$, and this contradicts $(2.1)$. Now, set
$$\delta^*=
\inf_{x\in B_r}\sup_{y\in B_r}\left({{\rho}\over {2}}\|x\|^2-{{1}\over {2}}\|x-y\|^2-P(x)\right)-
\sup_{y\in B_r}\inf_{x\in B_r}\left({{\rho}\over {2}}\|x\|^2-{{1}\over {2}}\|x-y\|^2-P(x)\right)\ .$$
In view of $(2.3)$, we have $\delta^*>0$.
Let us show that $\delta^*$ has the property described in the conclusion. So, fix any convex set $S\subset \hbox {\rm int}(B_r)$, dense
in $B_r$ and fix also any $C^1$ functional $\gamma:E\to {\bf R}$ satisfying 
$$\sup_{B_r}\gamma-\inf_{B_r}\gamma<\delta^*\ ,$$
$$\sup_{x\in B_r}\|\gamma'(x)\|\leq M\ ,$$
such that the restriction to $B_r$ of the functional $x\to \sigma\|x\|^2-P(x)-\gamma(x)$ is sequentially weakly lower semicontinuous.
We are going to apply Theorem 1.1 taking $X=Y=B_r$, $\alpha=0$, $\beta=-\gamma$ and
$$I(x,y)={{\rho}\over {2}}\|x\|^2\ ,$$
$$J(x,y)={{1}\over {2}}\|x-y\|^2+P(x)$$
for all $(x,y)\in E\times Y$. For what above, we have
$$\sup_Y\inf_X(I-J)<\inf_X\sup_Y(I-J)$$
and
$$\sup_X(\alpha+\beta)-\inf_X(\alpha+\beta)<\inf_X\sup_Y(I-J)-\sup_Y\inf_X(I-J)\ .$$
Moreover, conditions $(i_1), (i_3)$ and $(i_2)$ (taking as $\tau$ the weak topology on $X$)
are clearly satisfied. Finally, let us check that $(i_4)$ is satisfied. To this end, fix $u, y\in B_r$. We want to show that, for some $v\in B_r$,
we have
$$J'_x(u,y)+\gamma'(u)=I'(v)\ .$$
This means 
$$u-y+P'(u)+\gamma'(u)=\rho v$$
that is
$$\|u-y+P'(u)+\gamma'(u)\|\leq \rho r\ .\eqno{(2.5)}$$
Clearly
$$\|u-y+P'(u)+\gamma'(u)\|\leq \sup_{x\in B_r}\|x+P'(x)\|+r+M\ .\eqno{(2.6)}$$
From $(2.2)$, we get
$${{\sup_{x\in B_r}\|x+P'(x)\|+M}\over {r}}\leq {{1}\over {\lambda}}$$
and so
$${{\sup_{x\in B_r}\|x+P'(x)\|+r+M}\over {r}}\leq {{1}\over {\lambda}}+1=\rho$$
from which, together $(2.6)$, we get $(2.5)$. So, the assumptions of Theorem 1.1 are satisfied. Then, since $-S$ is a convex set dense in $B_r$,
there exists $\tilde z\in -S$, such that the equation 
$$(\rho-1)x+\tilde z-P'(x)-\gamma'(x)=0$$
has at least two solutions lying in $B_r$. So, the conclusion follows taking $\tilde y=-\tilde z$.\hfill $\bigtriangleup$\par
\medskip
Before presenting another application of Theorem 1.1, we prove a new general criterion ensuring the strict minimax inequality.\par
\medskip
THEOREM 2.2. - {\it Let $V$ be a topological space, let $H$ be a real Hilbert, let $Y$ be a closed ball in $H$ centered at $0$,
and let $Q:V\to {\bf R}$, $\psi:V\to H$. Assume that the functional $x\to Q(x)-\langle\psi(x),y\rangle$ is lower semicontinuous
for each $y\in Y$, while the functional $x\to Q(x)-\langle\psi(x),y_0\rangle$ is inf-compact for some $y_0\in Y$. Moreover,
assume that, for each $x\in V$, there exists $u\in V$ such that
$$Q(x)=Q(u)$$
and
$$\psi(x)=-\psi(u)\ .$$
Finally, assume that there is no global minimum of $Q$ at which $\psi$ vanishes.\par
Then, we have
$$\sup_{y\in Y}\inf_{x\in V}(Q(x)-\langle\psi(x),y\rangle)<\inf_{x\in V}\sup_{y\in Y}(Q(x)-\langle\psi(x),y\rangle)\ .$$}
\smallskip
PROOF. Arguing by contradiction, assume that
$$\sup_{y\in Y}\inf_{x\in V}(Q(x)-\langle\psi(x),y\rangle)=\inf_{x\in V}\sup_{y\in Y}(Q(x)-\langle\psi(x),y\rangle)\ .\eqno{(2.7)}$$
By Proposition 2.A (we consider $Y$ with the weak topology), $(2.7)$
would imply the existence of $(x^*,y^*)\in V\times Y$ such that
$$Q(x^*)-\langle\psi(x^*),y^*\rangle=\inf_{x\in V}(Q(x)-\langle\psi(x),y^*\rangle)=\sup_{y\in Y}(Q(x^*)-\langle\psi(x^*),y\rangle)\ .\eqno{(2.8)}$$
From this, we get
$$\langle\psi(x^*),y^*\rangle=\inf_{y\in Y}\langle\psi(x^*),y\rangle=-r\|\psi(x^*)\|$$
where $r$ is the radius of $Y$. Notice that $\psi(x^*)=0$. Indeed, if $\psi(x^*)\neq 0$, we would have
$$y^*=-{{r}\over {\|\psi(x^*)\|}}\psi(x^*)\ ,$$
and so $(2.8)$ would give
$$Q(x^*)+r\|\psi(x^*)\|=\inf_{x\in V}\left(Q(x)+{{r}\over {\|\psi(x^*)\|}}\langle\psi(x),\psi(x^*)\rangle\right)\ .\eqno{(2.9)}$$
Now, by assumption, there is a point $u^*\in V$ such that
$$Q(u^*)=Q(x^*)$$
and
$$\psi(u^*)=-\psi(x^*)\ .$$
Therefore, we have
$$Q(u^*)+{{r}\over {\|\psi(x^*)\|}}\langle\psi(u^*),\psi(x^*)\rangle=Q(x^*)-r\|\psi(x^*)\|<Q(x^*)+r\|\psi(x^*)\|\ ,$$
contradicting $(2.9)$. Hence, $\psi(x^*)=0$. But then, $(2.8)$ gives
$$Q(x^*)=\inf_{x\in V}(Q(x)-\langle\psi(x),y^*\rangle)\ .\eqno{(2.10)}$$
By assumption, $x^*$ is not a global minimum of $Q$. So, there exists $x_1\in V$ such that $Q(x_1)<Q(x^*)$. By assumption again, 
there exists $u_1\in V$ such that
$$Q(u_1)=Q(x_1)$$
and
$$\psi(u_1)=-\psi(x_1)\ .$$
Clearly, we have either
$$Q(x_1)-\langle\psi(x_1),y^*\rangle\leq Q(x_1)$$
or
$$Q(u_1)-\langle\psi(u_1),y^*\rangle\leq Q(x_1)\ .$$
Hence, we would have
$$\inf_{x\in V}(Q(x)-\langle\psi(x),y^*\rangle)<Q(x^*)\ ,$$
contradicting $(2.10)$.\hfill $\bigtriangleup$
\medskip
THEOREM 2.3. -  {\it Let $E, H$ be two real Hilbert spaces, let $r, s>0$, let $P:E\to {\bf R}$ be a $C^1$ 
even functional such that the restriction to $B_r$ of the functional $x\to \sigma\|x\|^2-P(x)$ is sequentially weakly lower semicontinuous for some 
$\sigma>0$, and let $\psi:E\to H$ be a sequentially weakly continuous $C^1$ odd operator.\par
Then, for each $M\geq 0$ and for each $\lambda$ satisfying
$$0<\lambda\leq \min\left\{{{r}\over {\sup_{x\in B_r}(\|P'(x)\|
+s\|\psi'(x)\|_{{\cal L}(E,H)})+M}}, {{1}\over {2\sigma}}\right\}$$
and such that the restriction to $B_r$ of the functional $x\to {{1}\over {2\lambda}}\|x\|^2-P(x)$ has no global minimum at which $\psi$ vanishes,
there exists $\delta^*>0$ with the following property: 
for each convex set $S\subseteq H$, dense in $H$, and for each
$C^1$ functional $\gamma:E\to {\bf R}$ satisfying 
$$\sup_{B_r}\gamma-\inf_{B_r}\gamma<\delta^*\ ,$$
$$\sup_{x\in B_r}\|\gamma'(x)\|\leq M\ ,$$
such that the restriction to $B_r$ of the functional $x\to \sigma\|x\|^2-P(x)-\gamma(x)$ is sequentially weakly lower semicontinuous,
there exists $\tilde y$ in $S$, with $\|\tilde y\|_H<s$, such that the equation
$$x=\lambda\left(P'(x)+\gamma'(x)+{{d}\over {dx}}\langle \psi(x),\tilde y\rangle_H\right)$$
has at least two solutions lying in $B_r$.}\par
\smallskip
PROOF. Fix $M\geq0$ and fix $\lambda$ such that
$$0<\lambda\leq \min\left\{{{r}\over {\sup_{x\in B_r}(\|P'(x)\|
+s\|\psi'(x)\|_{{\cal L}(E,H)})+M}}, {{1}\over {2\sigma}}\right\}\eqno{(2.11)}$$
and that the restriction to $B_r$ of the functional $x\to {{1}\over {2\lambda}}\|x\|^2-P(x)$ has no global minimum at which $\psi$ vanishes.
Denote by $Y$ the closed ball in $H$, centered at $0$, of radius $s$. Since the operator $\psi$ is odd and
the functional $x\to {{1}\over {2\lambda}}\|x\|^2-P(x)$ is even and
its restriction to $B_r$ is sequentially weakly lower semicontinuous, we can apply Theorem 2.2, taking $V=B_r$ with the weak topology. As a consequence, the number
$$\delta^*:=\inf_{x\in B_r}\sup_{y\in Y}\left({{1}\over {2\lambda}}\|x\|^2-P(x)-\langle\psi(x),y\rangle_H\right)-\sup_{y\in Y}\inf_{x\in B_r}
\left({{1}\over {2\lambda}}\|x\|^2-P(x)-\langle\psi(x),y\rangle_H\right)$$
is positive. Let us show that $\delta^*$ has the property described in the conclusion. So, fix any convex set $S\subset Y$, dense
in $Y$ and fix also any $C^1$ functional $\gamma:E\to {\bf R}$ satisfying 
$$\sup_{B_r}\gamma-\inf_{B_r}\gamma<\delta^*\ ,$$
$$\sup_{x\in B_r}\|\gamma'(x)\|\leq M\ ,$$
such that the restriction to $B_r$ of the functional $x\to \sigma\|x\|^2-P(x)-\gamma(x)$ is sequentially weakly lower semicontinuous.
We are going to apply Theorem 1.1 taking $X=B_r$, $\alpha=0$, $\beta=-\gamma$ and
$$I(x,y)={{1}\over {2\lambda}}\|x\|^2\ ,$$
$$J(x,y)=\langle\psi(x),y\rangle_H+P(x)$$
for all $(x,y)\in E\times Y$. For what above, we have
$$\sup_Y\inf_X(I-J)<\inf_X\sup_Y(I-J)$$
and
$$\sup_X(\alpha+\beta)-\inf_X(\alpha+\beta)<\inf_X\sup_Y(I-J)-\sup_Y\inf_X(I-J)\ .$$
Moreover, conditions $(i_1), (i_3)$ and $(i_2)$ (taking as $\tau$ the weak topology on $X$)
are clearly satisfied. Finally, let us check that $(i_4)$ is satisfied. To this end, fix $u\in B_r, y\in Y$. We want to show that, for some $v\in B_r$,
we have
$$J'_x(u,y)+\gamma'(u)=I'(v)\ .$$
This means 
$$\sup_{\|w\|\leq 1}|\langle\psi'(u)(w),y\rangle_H+\langle P'(u),w\rangle+\langle\gamma'(u),w\rangle|\leq {{r}\over {\lambda}}\ .\eqno{(2.12)}$$
Clearly
$$\sup_{\|w\|\leq 1}|\langle\psi'(u)(w),y\rangle_H+\langle P'(u),w\rangle+\langle\gamma'(u),w\rangle|\leq
s\|\psi'(u)\|_{{\cal L}(E,H)}+\|P'(u)\|+\|\gamma'(u)\|$$
$$\leq \sup_{x\in B_r}(\|P'(x)\|+s\|\psi'(x)\|_{{\cal L}(E,H)})+M$$
and so $(2.12)$ follows from $(2.11)$.
Therefore, the assumptions of Theorem 1.1 are satisfied. Then, since $-S\cap \hbox {\rm int}(Y)$ is a convex set dense in $Y$,
there exists $\tilde z\in -S$, such that the equation 
$${{1}\over {\lambda}}x-{{d}\over {dx}}\langle\psi(x),\tilde y\rangle_H-P'(x)-\gamma'(x)=0$$
has at least two solutions lying in $B_r$. So, the conclusion follows taking $\tilde y=-\tilde z$.\hfill $\bigtriangleup$\par
\medskip
We now present two applications of Theorems 2.1 and 2.3 to elliptic equations involving nonlinearities with critical growth.\par
\smallskip
In the sequel, $\Omega\subset {\bf R}^n$ ($n\geq 3$) is a smooth bounded domain.\par
\smallskip
We denote by ${\cal A}$ (resp. $\tilde {\cal A}$) the class of all
Carath\'eodory functions $f:\Omega\times {\bf R}\to {\bf R}$ such that
$$\sup_{(x,\xi)\in \Omega\times {\bf R}}{{|f(x,\xi)|}\over
{1+|\xi|^q}}<+\infty\ ,$$
where  $0<q\leq {{n+2}\over {n-2}}$ (resp. $0<q<{{n+2}\over {n-2}}$ ).
\smallskip
We consider the Sobolev space $H^1_0(\Omega)$ endowed with the scalar product
$$\langle u,v\rangle=\int_{\Omega}\nabla u(x)\nabla v(x)dx\ .$$
\smallskip
As usual, we denote by $H^{-1}(\Omega)$ the dual space of $H^1_0(\Omega)$.\par
\smallskip
Given $f\in {\cal A}$ and $\varphi\in H^{-1}(\Omega)$,  consider the following Dirichlet problem
$$\cases {-\Delta u= f(x,u) +\varphi
 & in
$\Omega$\cr & \cr u=0 & on
$\partial \Omega$\ .\cr}\eqno{(P_{f,\varphi})} $$
 Let us recall
that a weak solution
of $(P_{f,\varphi})$ is any $u\in H^1_0(\Omega)$ such that
 $$\int_{\Omega}\nabla u(x)\nabla v(x)dx
-\int_{\Omega}f(x,u(x))v(x)dx-\varphi(v)=0$$
for all $v\in H^1_0(\Omega)$.\par
\smallskip
By classical results, if we set $F(x,\xi)=\int_0^{\xi}f(x,t)dt$, the functional $u\to \int_{\Omega}F(x,u(x))dx$ is $C^1$ in $H^1_0(\Omega)$
and 
the weak solutions of problem $(P_{f,\varphi})$ agree with the critical points in $H^1_0(\Omega)$ of the functional
$$u\to {{1}\over {2}}\int_{\Omega}|\nabla u(x)|^2dx-\int_{\Omega}F(x,u(x))dx-\varphi(u)\ .$$
Recall also that the functional $u\to \int_{\Omega}F(x,u(x))dx$ is sequentially weakly continuous in $H^1_0(\Omega)$ provided that $f\in
\tilde {\cal A}$. This is not the case for certain $f\in {\cal A}$.

Here is the first application of Theorem 2.1.\par
\medskip
THEOREM 2.4. - {\it Let $f, g\in \tilde {\cal A}$. Assume that, for each $x\in \Omega$, the function $f(x,\cdot)$ is odd and that
$$\lim_{\xi\to 0^+}{{\inf_{x\in \Omega}F(x,\xi)}\over {\xi^2}}=+\infty\ .\eqno{(2.13)}$$
Then, for each $r>0$, there exists $\lambda^*>0$ with the following property: for each $\lambda\in ]0,\lambda^*[$, for each $\mu\in ]-\lambda^*,
\lambda^*[$ and for each convex set $S\subseteq H^{-1}(\Omega)$, dense in $H^{-1}(\Omega)$, there exists $\tilde\varphi\in S$,
with $\|\tilde\varphi\|_{H^{-1}(\Omega)}<r$, such that the problem
$$\cases{-\Delta u=\lambda(|u|^{{{4}\over {n-2}}}u+f(x,u)+\mu g(x,u) +\tilde\varphi)
& in
$\Omega$\cr & \cr u=0 & on
$\partial \Omega$ \cr}$$
has at least two weak solutions whose norms in $H^1_0(\Omega)$ are less than or equal to $r$.}\par
\smallskip
PROOF. We are going to apply Theorem 2.1 taking $E=H^1_0(\Omega)$ and 
$$P(u)={{n-2}\over {2n}}\int_{\Omega}|u(x)|^{{2n}\over {n-2}}dx+\int_{\Omega}F(x,u(x))dx\ .$$
Clearly, $P$ is even since $F(x,\cdot)$ is so. Moreover, $(2.1)$ follows easily from $(2.13)$ and the fact that, for some $\sigma>0$,
the restriction to $B_r$ of the functional $u\to \sigma\|u\|^2-P(u)$ is sequentially weakly lower semicontinuous is guaranteed by Lemma 1 of [1].
Next, consider the functional $\gamma:H^1_0(\Omega)\to {\bf R}$ defined by
$$\gamma(u)=\int_{\Omega}G(x,u(x))dx\ ,$$
where $G(x,\xi)=\int_0^{\xi}g(x,t)dt$. Apply Theorem 2.1 with $M=\sup_{u\in B_r}\|\gamma'(u)\|$. Accordingly, there exists $\eta_r>0$ with
the following property: for each $\lambda\in ]0,\eta_r]$ there exists $\delta^*>0$ such that, for each $\mu\in {\bf R}$,
with $|\mu|<{{\delta^*}\over {\sup_{B_r}\gamma-\inf_{B_r}\gamma}}$, and for each convex set $S\subseteq H^{-1}(\Omega)$, dense in $H^{-1}(\Omega)$, there exists $\tilde\varphi\in S$, with $\|\tilde\varphi\|_{H^{-1}(\Omega)}<r$, such that the functional
$$u\to {{1}\over {2}}\int_{\Omega}|\nabla u(x)|^2-\lambda\left({{n-2}\over {2n}}\int_{\Omega}|u(x)|^{{2n}\over {n-2}}dx+\int_{\Omega}F(x,u(x))dx
-\mu\int_{\Omega}G(x,u(x))dx-\tilde \varphi(u)\right)$$
has at least two critical points in $H^1_0(\Omega)$ lying in $B_r$. So, we can take $\lambda^*=\min\left\{\eta_r,{{\delta^*}\over {\sup_{B_r}\gamma-\inf_{B_r}\gamma}}\right\}$.\hfill $\bigtriangleup$\par
\medskip
REMARK 2.1. - We refer to the very recent [3] (and references therein) for other remarkable results on nonhomogeneous elliptic equations
involving nonlinearities with critical growth.\par
\medskip
The second application of Theorem 2.1 is as follows.\par
\medskip
THEOREM 2.5. - {\it Let $f, g, h\in \tilde {\cal A}$. Assume that, for each $x\in \Omega$, the function $f(x,\cdot)$ is odd and the function
$h(x,\cdot)$ is even and vanishes at most at $0$. Also, suppose that
$$\lim_{\xi\to 0^+}{{\inf_{x\in \Omega}F(x,\xi)}\over {\xi^2}}=+\infty\ .$$
Then, for each $r, s>0$, there exists $\lambda^*>0$ with the following property: for each $\lambda\in ]0,\lambda^*[$, for each $\mu\in ]-\lambda^*,
\lambda^*[$ and for each convex set $S\subseteq L^{\infty}(\Omega)$, dense in $L^2(\Omega)$, there exists $\tilde\varphi\in S$,
with $\|\tilde\varphi\|_{L^2(\Omega)}<s$, such that the problem
$$\cases{-\Delta u=\lambda(|u|^{{{4}\over {n-2}}}u+f(x,u)+\mu g(x,u) +\tilde\varphi(x)h(x,u))
& in
$\Omega$\cr & \cr u=0 & on
$\partial \Omega$ \cr}$$
has at least two weak solutions whose norms in $H^1_0(\Omega)$ are less than or equal to $r$.}\par
\smallskip
PROOF. We are going to apply Theorem 2.1 taking $E=H^1_0(\Omega)$, $H=L^2(\Omega)$, 
$$P(u)={{n-2}\over {2n}}\int_{\Omega}|u(x)|^{{2n}\over {n-2}}dx+\int_{\Omega}F(x,u(x))dx$$
and
$$\psi(u)(x)=H(x,u(x))$$
for all $u\in H^1_0(\Omega)$, $x\in \Omega$, where $H(x,\xi)=\int_0^{\xi}h(x,t)dt.$
We already noticed in the proof of Theorem 2.4 that $P$ satisfies the requirements of Theorem 2.1. Concerning $\psi$, it is odd (since
$h(x,\cdot)$ is even) and sequentially weakly continuous by the Rellich-Kondrachov theorem. Fix  $r, s>0$.
Notice that $\psi(u)=0$ if and only
if $u=0$ and, for each $\nu>0$, we 
$$\inf_{u\in B_r}\left(\nu\int_{\Omega}|\nabla u(x)|^2dx-P(u)\right)<0$$
since
$$\lim_{u\to 0}{{P(u)}\over {\int_{\Omega}|\nabla u(x)|^2dx}}=+\infty\ .$$
Consequently, there is no global minimum of the restriction to $B_r$ of the functional $u\to \nu\int_{\Omega}|\nabla u(x)|^2dx-P(u)$
at which $\psi$ vanishes. Now, let $\gamma$ be as in the proof of Theorem 2.4 and apply Theorem 2.1 with 
$M=\sup_{u\in B_r}\|\gamma'(u)\|$. 
Accordingly, there exists $\eta_r>0$ with
the following property: for each $\lambda\in ]0,\eta_r]$ there exists $\delta^*>0$ such that, for each $\mu\in {\bf R}$,
with $|\mu|<{{\delta^*}\over {\sup_{B_r}\gamma-\inf_{B_r}\gamma}}$, and for each convex set $S\subseteq L^{\infty}(\Omega)$, dense in 
in $L^2(\Omega)$, there exists $\tilde\varphi\in S$, with $\|\tilde\varphi\|_{L^2(\Omega)}<r$, such that the functional
$$u\to {{1}\over {2}}\int_{\Omega}|\nabla u(x)|^2-\lambda\left({{n-2}\over {2n}}\int_{\Omega}|u(x)|^{{2n}\over {n-2}}dx+\int_{\Omega}F(x,u(x))dx
-\mu\int_{\Omega}G(x,u(x))dx-\int_{\Omega}\tilde\varphi(x)H(x,u(x))\right)$$
has at least two critical points in $H^1_0(\Omega)$ lying in $B_r$. So, we can take $\lambda^*=\min\left\{\eta_r,{{\delta^*}\over {\sup_{B_r}\gamma-\inf_{B_r}\gamma}}\right\}$.\hfill $\bigtriangleup$\par
\bigskip
\bigskip
\bigskip
{\bf Acknowledgements:} This work has been funded by the European Union - NextGenerationEU Mission 4 - Component 2 - Investment 1.1 under the Italian Ministry of University and Research (MUR) programme "PRIN 2022" - grant number 2022BCFHN2 - Advanced theoretical aspects in PDEs and their applications - CUP: E53D23005650006. The author has also been supported by the Gruppo Nazionale per l'Analisi Matematica, la Probabilit\`a e 
le loro Applicazioni (GNAMPA) of the Istituto Nazionale di Alta Matematica (INdAM) and by the Universit\`a degli Studi di Catania, PIACERI 2020-2022, Linea di intervento 2, Progetto ”MAFANE”.
\bigskip
\bigskip
\bigskip
\bigskip
\centerline {\bf References}\par
\bigskip
\bigskip
\noindent
[1]\hskip 5pt F. FARACI and C. FARKAS, {\it A quasilinear elliptic problem involving critical Sobolev
exponents}, Collect. Math., {\bf 66} (2015), 243-259.\par
\smallskip
\noindent
[2]\hskip 5pt A. MOAMENI, {\it Critical point theory on convex subsets with applications in differential equations and analysis},
 J. Math. Pures Appl., {\bf 141} (2020), 266-315.\par
\smallskip
\noindent
[3]\hskip 5pt K. PERERA, {\it A general perturbation theorem with applications to nonhomogeneous critical growth elliptic problems},
J. Differential Equations, {\bf 389} (2024), 150-189.\par
\smallskip
\noindent
[4]\hskip 5pt B. RICCERI, {\it On a minimax theorem: an improvement, a new proof and an overview of its applications},
Minimax Theory Appl., {\bf 2} (2017), 99-152.\par
\smallskip
\noindent
[5] \hskip 5pt B. RICCERI, {\it A more complete version of a minimax theorem}, Appl. Anal. Optim., {\bf 5} (2021), 251-261.\par
\smallskip
\noindent
[6]\hskip 5pt B. RICCERI, {\it Addendum to  ``A more complete version of a minimax theorem"}, Appl. Anal. Optim., {\bf 6} (2022), 195-197.\par
\smallskip
\bigskip
\bigskip
\bigskip
\bigskip
Department of Mathematics and Informatics\par
University of Catania\par
Viale A. Doria 6\par
95125 Catania, Italy\par
{\it e-mail address}: ricceri@dmi.unict.it

\bye